\documentclass[11pt]{amsart}
\usepackage{fullpage,amssymb,amsmath,verbatim,amstext,eucal}

\newtheorem{lemma}{Lemma}[section] 

\newtheorem{theorem}[lemma]{Theorem}
\newtheorem{corollary}[lemma]{Corollary}

\newenvironment{remark}[1]{\refstepcounter{lemma}%
\vskip 5pt \par\noindent {\bf #1\ \thelemma .}}{\vskip 5pt \par}

\newenvironment{remark*}[1]{\par \vskip 5pt \noindent 
{\bf #1.}}{\vskip 5pt \par}

\newlength{\dqlength}

\newcommand{\Ad}{\operatorname{Ad}}

\newcommand{\bh}{\ensuremath{{\mathcal B}({\mathcal H})}}

\newcommand{\condnum}[1]{\norm{{#1}^{\vphantom{1}}}\norm{\inv{#1}}}
\newcommand{\cstar}{\hbox{$C^*$}}
\newcommand{\cstaralg}{$C^*$-algebra}
\newcommand{\cstarenv}{\hbox{$C^*_{\rm env}$}}

\newcommand{\dist}{\operatorname{dist}}

\newcommand{\ds}{\displaystyle}
\newcommand{\eps}{\mbox{$\varepsilon$}}

\newcommand{\inv}[1]{{#1}^{-1}}

\newcommand{\lat}{{\rm Lat}\,}

\newcommand{\norm}[1]{\left\|{#1}\right\|}
\newcommand{\normcb}[1]{\left\|{#1}\right\|_{\rm{cb}}}

\providecommand{\qed}%
{\hfill \vrule height5pt width4pt depth1pt \vspace{+2.00ex}}

\newcommand{\bbA}{{\mathbb{A}}}

\newcommand{\bbC}{{\mathbb{C}}}
\newcommand{\bbD}{{\mathbb{D}}}

\newcommand{\bbN}{{\mathbb{N}}}

\newcommand{\bbT}{{\mathbb{T}}}

  \newcommand{\A}{{\mathcal{A}}}
  
  \newcommand{\B}{{\mathcal{B}}}
  \newcommand{\C}{{\mathcal{C}}}
  
  \newcommand{\D}{{\mathcal{D}}}

\renewcommand{\H}{{\mathcal{H}}}
  
  \newcommand{\J}{{\mathcal{J}}}
  \newcommand{\K}{{\mathcal{K}}}  

  \newcommand{\M}{{\mathcal{M}}}

\renewcommand{\P}{{\mathcal{P}}}

\newcommand{\fM}{{\mathfrak{M}}}

\begin{document}

\title{Norming Algebras and Automatic Complete
Boundedness of Isomorphisms of Operator Algebras}

\author[D.R. Pitts]{David R. Pitts}
\address{Dept. of Mathematics\\
University of Nebraska-Lincoln\\ Lincoln, NE\\ 68588-0130}
\email{dpitts2@math.unl.edu}
\subjclass[2000]{47L30, 46L07, 47L55}
\keywords{completely bounded isomorphism, \cstaralg, operator algebra}
\date{September 20, 2006}

\begin{abstract}
We combine the notion of norming algebra introduced by Pop, Sinclair
and Smith with a result of Pisier to show that if $\A_1$ and $\A_2$
are operator algebras, then any bounded epimorphism of $\A_1$ onto
$\A_2$ is completely bounded provided that $\A_2$ contains a norming
\cstar-subalgebra.  We use this result to give some insights into
Kadison's Similarity Problem: we show that every faithful bounded
homomorphism of a \cstaralg\ on a Hilbert space has completely bounded
inverse, and show that a bounded representation of a \cstaralg\ is
similar to a $*$-representation precisely when the image operator
algebra $\lambda$-norms itself.  We give two applications to isometric
isomorphisms of certain operator algebras.  The first is an extension
of a result of Davidson and Power on isometric isomorphisms of CSL
algebras.  Secondly, we show that an isometric isomorphism between
subalgebras $\A_i$ of \cstar-diagonals $(\C_i,\D_i)$ ($i=1,2$)
satisfying $\D_i\subseteq\A_i\subseteq \C_i$ extends uniquely to a
$*$-isomorphism of the \cstaralg s generated by $\A_1$ and $\A_2$;
this generalizes results of Muhly-Qiu-Solel and Donsig-Pitts.
 
\end{abstract}
\maketitle

\vskip 12pt

\section{Introduction and Norming Algebras}
Let $\A$ be a unital operator algebra and $u:\A\rightarrow \bh$ be a
homomorphism.  If $u$ is contractive (resp.\ isometric or bounded), it
is not generally the case that $u$ is completely contractive (resp.\
completely isometric or completely bounded).  However, in some cases
it is possible to conclude that if $u$ is isometric or contractive,
then $u$ is completely isometric or completely contractive.  For example,
the contractive homomorphisms of a \cstaralg\ are exactly the $*$
homomorphisms, and hence any contractive homomorphism of a \cstaralg\
into \bh\ is completely contractive.  It is not known however if every
bounded representation of a \cstaralg\ is completely bounded; a result
of Haagerup (stated as Theorem~\ref{HaagerupSimThm} below)
shows this question is equivalent to Kadison's similarity problem.

The purpose of this note is to give a sufficient condition on an
operator algebra $\B$ which ensures that every bounded epimorphism
$u:\A\rightarrow \B$ of the operator algebra $\A$ onto $\B$ is
completely bounded.  We show that this condition can be used to give
simple proofs of several results in the literature, and also that
every bounded faithful representation of a \cstaralg\ is  bounded below.


Throughout, all operator algebras are norm-closed, and we typically
use  $\A$ and $\B$ to denote operator algebras.
(There are several texts containing the background we need
for operator spaces and operator algebras, see for
example~\cite{BlecherLeMerdyOpAlThMo,EffrosRuanOpSp,PaulsenCoBoMaOpAl}.)
An \textit{operator
$\A$-$\B$-bimodule} is an operator space $\M$ which is a left-$\A$,
right-$\B$ bimodule where the bimodule action is completely
contractive in the sense that
for any $A\in
M_{np}(\A)$, $M\in M_p(\M)$ and $B\in M_{pn}(\B)$,
\[ \norm{AMB}_{M_n(\M)}\leq
\norm{A}_{M_{np}(\A)}\norm{M}_{M_p(\M)}\norm{B}_{M_{pn}(\B)}.\] (We
shall generally not write subscripts on the norms in the sequel,
unless necessary for clarity.)  We will sometimes write $_A\M_\B$ when
$\M$ is 
an $A$-$\B$-bimodule.

For convenience, we will always assume that
operator algebras are unital unless explicitly stated otherwise.
However, many of the results in the sequel are valid for non-unital
algebras. 

Given such an $\A$-$\B$ operator bimodule $\M$, we may define a family
of norms $\eta_n$ on $M_n(\M)$ as follows:  for $X\in M_n(\M)$,
\[\eta_n(X):=\sup\{\norm{RXC}: R\in M_{1n}(\A), C\in M_{n1}(\B) \text{
  and } \max\{\norm{C},\norm{R}\}\leq 1\}.\]  Clearly $\eta_n(X)\leq
  \norm{X}_{\M_n(\M)}.$

\begin{remark}{Definition}(\cite{PopSinclairSmithNoC*Al})\label{normingdef} 
 Let $\lambda>0$ be a real number.  We say that $_\A\M_\B$ is
  \textit{$\lambda$-normed} by $\A$ and $\B$ if for every $n\in \bbN$
  and $X\in M_n(\M)$,
\[\lambda\norm{X}\leq \eta_n(X).\]
When  $\lambda=1$,
we say that $_\A\M_\B$ is \textit{normed} by $\A$-$\B$.  When $\A=\B$,
we simply say that $\M$ is normed by $\A$.
\end{remark}  

When $\A$ and $\B$ are \cstaralg s, C. Pop~\cite{PopBiNoReEsHi} shows
that $(\M,\{\eta_n\})$ is the smallest operator space structure on
the Banach space $\M$ which is compatible
with the module structure of $\M$.

We will make essential use of the following result.

\begin{theorem}[{\cite[Lemma~7.7, p.\ 128]{PisierSiPrCoBoMaIIed}}] 
\label{PisierLemma} 
Let $\A$ be a
  \cstaralg\ and suppose that $u:\A\rightarrow \bh$ is a bounded
  homomorphism.  Then for any $R\in M_{1,n}(\A)$ and $C\in
  M_{n,1}(\A)$ we have
\begin{align*}
\norm{u_{1,n}(R)(u_{1,n}(R))^*}&\leq
\norm{u}^4\norm{RR^*}\qquad\text{and}\\
\norm{(u_{n,1}(C))^*u_{n,1}(C)}&\leq \norm{u}^4\norm{C^*C}.
\end{align*}
\end{theorem}
\begin{remark}{Remark}
 Actually,~\cite[Lemma~7.7, p.\ 128]{PisierSiPrCoBoMaIIed} gives the
 statement and proof for the column $C$.  However to obtain the
 statement for $R$, one applies the statement for $C$ to the opposite
 algebra $\A^{\rm op}$ and the homomorphism $\psi$ of $\A^{\rm op}$ on
 the conjugate Hilbert space $\overline{\H}$ given by $d\in\A\mapsto
 \psi(d)$, where $\psi(d)\overline{h}=\overline{(u(d))^*h}$.
\end{remark}

The following is our central observation.  The key idea is to
combine the techniques of~\cite[Theorem~2.10]{PopSinclairSmithNoC*Al}
with Theorem~\ref{PisierLemma}.  (See
also~\cite[Theorem~2.1]{SmithCoBoMoMaHaTePr}.)

\begin{theorem}\label{autocb}   Let $\B$ be a norm-closed operator 
  algebra which contains a \cstaralg\ $\D$,  and assume $\D$ 
  norms $\B$.  If $\A$ is a norm-closed operator
  algebra and $u:\A\rightarrow \B$ is a bounded isomorphism,
  then $u$ is completely bounded and 
$$\norm{u}_{cb}\leq \norm{u}\norm{u^{-1}}^4.$$
\end{theorem}

\begin{remark*}{Remarks}  If $\D$ only $\lambda$-norms $\B$, then
  $u$ is completely bounded with $\norm{u}_{cb}\leq
  \lambda^{-1}\norm{u}\norm{u^{-1}}^4$.   Also, clearly the theorem
  applies when $u$ is a bounded epimorphism:  simply  replace $\A$ with
  $\A/\ker u$ and $u$ with the induced isomorphism of the quotient
  onto $\B$. 
\end{remark*}

\begin{proof} Let $T=(T_{ij})\in M_n(\A)$.  Then for any 
$R\in M_{1n}(\D)$ and $C\in
M_{n1}(\D)$ with $\norm{R}\leq 1$ and $\norm{C}\leq 1$ we have, 
\begin{align*}
\norm{Ru_n(T)C} &= \norm{\sum_{i,j=1}^nR_i u(T_{ij})C_j}\\
&\leq\norm{u}\norm{\sum_{i,j=1}^n
   u^{-1}(R_i)T_{ij} u^{-1}(C_j)}\\
&=\norm{ u}\norm{ u^{-1}_{1,n}(R) T  u^{-1}_{n,1}(C)}\\
&\leq \norm{ u}  \norm{ u^{-1}_{1,n}(R)} \norm{T} 
\norm{ u^{-1}_{n,1}(C)}.
\end{align*}
 We may assume that $\A$ is represented completely isometrically as
operators acting on a Hilbert space $\H$.  Thus, $ u^{-1}|_\D$ is
a bounded homomorphism of $\D$ into $\bh$.  By
Lemma~\ref{PisierLemma}, 
$\norm{ u^{-1}_{n,1}(C)} \leq
\norm{ u^{-1}}^2$ and 
$\norm{ u^{-1}_{1,n}(R)}\leq \norm{ u^{-1}}^2$.  Taking
suprema over $R$ and $C$ gives
$$\norm{ u_n(T)}\leq \norm{ u}\norm{ u^{-1}}^4 \norm{T},$$
as desired.  
\end{proof}

When the isomorphism is isometric, more can be said.  For any operator
algebra $\B$, let $\cstarenv(B)$ be the $\cstar$-envelope of $\B$.

\begin{corollary}\label{isometriciso}  For $i=1,2$, suppose that
  $\A_i$ are operator algebras and $\D_i\subseteq \C_i$ is a
  norming \cstar-subalgebra of $\A_i$.  If
  $ u:\A_1\rightarrow \A_2$ is an isometric isomorphism, then
  $ u$ extends uniquely to a $*$-isomorphism $\tilde{u}:
  \cstarenv(\A_1)\rightarrow \cstarenv(\A_2).$

\end{corollary}

\begin{proof}
Theorem~\ref{autocb} shows that $ u$ and $ u^{-1}$ are
complete contractions, so that $ u$ is a complete isometry.  The
result follows from the universal property of \cstar-envelopes.
\end{proof}

\section{Applications} 
In this section we record some consequences of Theorem~\ref{autocb}.  
We shall require the following closely related results of Haagerup and Paulsen.

\begin{theorem}[Haagerup~\cite{HaagerupSoSiPrCyRe}]
  \label{HaagerupSimThm}
 Suppose $\A$ is a
  \cstaralg\ and $ u:\A\rightarrow \bh$ is a completely bounded
  homomorphism.  Then there exists an invertible operator $S\in\bh$
  with $\condnum{S}=\normcb{ u}$ such that for every $a\in\A$, 
\[a\mapsto S u(a)S^{-1}\] is completely contractive (and hence a
  $*$-representation).
\end{theorem}

\begin{theorem}[Paulsen~\cite{PaulsenCoBoHoOpAl}] 
 \label{PaulsenSimThm}  Suppose $\A$ is a
  unital operator algebra and $ u:\A\rightarrow \bh$ is a
  completely bounded unital homomorphism.  Then there exists an invertible
  operator $S\in\bh$ with $\condnum{S}=\normcb{ u}$ such that for
  every $a\in\A$,
\[a\mapsto S u(a)S^{-1}\] is a completely contractive
  homomorphism. 
\end{theorem}

\subsection{Applications to \cstaralg s and Kadison's
  Similarity Problem}

 We begin with a new proof of a result
 of Gardner.

\begin{theorem}[Gardner~\cite{GardnerOnIsC*Al}]  Suppose $\A$ and $\B$
  are \cstaralg s and $ u:\A\rightarrow \B$ is an isomorphism.
  Then $ u$ is (completely) bounded and there exists a $*$-isomorphism
  $\alpha:\A\rightarrow\B$ and a bounded automorphism $\beta$ of $\B$
  such that $ u=\beta\circ\alpha$.  If $\B\subseteq\bh$, there
  exists a positive invertible operator $S\in\bh$ with
  $\condnum{S}\leq \norm{ u}$ so that $\beta=\Ad
  S$.
\end{theorem}
\begin{proof}
By~\cite[Lemma~2.3(i)]{PopSinclairSmithNoC*Al}, $\B$ norms itself, and
as observed by Gardner, $ u$ is bounded.  (One can also use a
result of B. Johnson~\cite{JohnsonUnCoNoTo} (see
also~\cite{SinclairAuCoLiOp})  concerning automatic continuity of a
homomorphism from a Banach algebra onto a semi-simple Banach algebra.)
Theorem~\ref{autocb} implies that $ u$ is completely bounded.
Assume that $\B\subseteq \bh$.  By Theorem~\ref{HaagerupSimThm}, there
exists an invertible operator $T$ with $\condnum{T}=\normcb{ u}$
such that $\Ad T\circ u$ is a completely contractive homomorphism,
so that $\Ad T\circ  u$ is therefore a $*$-homomorphism.  Let
$S=|T|,$ and let $U$ be the polar part of $T$, so $T=US$.  Then $(\Ad
S)(\B)$ is a \cstaralg, and it follows that $S^2\B S^{-2}=\B$.  By
Gardner's Invariance Theorem~\cite[Theorem~3.5]{GardnerOnIsC*Al},
$\beta:=\Ad S^{-1}$ is an automorphism of $\B$.  Now let $\alpha=\Ad
(U^*T)\circ  u$.  Since $\Ad T\circ  u$ is a $*$-isomorphism
of $\A$ onto its range, $\alpha$ is a $*$-isomorphism of $\A$ onto
$\B$.  Then $ u=\beta\circ\alpha.$
\end{proof}

\begin{remark}{Remark}  Gardner's original arguments give somewhat more.  In
  particular, he shows that if $\A$ and $\B$ are faithfully
  represented using the universal atomic representations, then
  $\alpha$ can be taken to have the form $\Ad U$ for some unitary $U$.
\end{remark}

The following is an immediate corollary of Theorem~\ref{autocb}.
\begin{theorem}\label{isocstar}
  Suppose $\A$ is an operator algebra and $\B$ is a
  \cstaralg.  If $ u:\A\rightarrow\B$ is an isomorphism, then
  $ u$ is automatically completely bounded and
  $\norm{ u}_{cb}\leq \norm{ u}\norm{ u^{-1}}^4$.
\end{theorem}
\begin{proof}
Continuity of $ u$ follows from Johnson's
theorem~\cite{JohnsonUnCoNoTo}, then complete boundedness follows from
Theorem~\ref{autocb} together with the fact that a \cstaralg\ norms
itself~\cite[Lemma~2.3(i)]{PopSinclairSmithNoC*Al}.
\end{proof}

We now wish to make some observations regarding Kadison's Similarity
Problem.  Recall that this problem asks whether every bounded
representation of a \cstaralg\ is similar to a $*$-representation,
which as noted above, is equivalent to the question of whether bounded
representations of \cstaralg s are automatically completely bounded.
Theorem~\ref{isocstar} can be used to prove
that bounded representations of \cstaralg s are (modulo the kernel)
``completely bounded below.''
\begin{theorem}\label{isocor} Suppose $\A$ is a \cstaralg\ and
  $u:\A\rightarrow \bh$ is a bounded homomorphism.  Put $\J=\ker u$.
 Then there exists a real number $k>0$ such that for every $n\in\bbN$
 and every $T\in M_n(\A)$,
$$k\dist(T,M_n(\J)) \leq \norm{u_n(T)}.$$
\end{theorem}
\begin{proof}
Since $\A/\J$ is a \cstaralg\ isomorphic under the map induced by $u$
to $u(\A)$, without loss of generality we may assume that $u$ is
one-to-one.  Thus our task is to prove that $ u:=u^{-1}$ is
completely bounded.  This will follow from Theorem~\ref{isocstar} once
we prove that the image $\B:=u(\A)$ is closed, and hence an operator
algebra.  We may also assume that $\A$ is unital and that $u(I)=I$.

Let $x$ be a non-zero element of $\A$, and let $\C$ be the unital
\cstar-subalgebra of $\A$ generated by $x^*x$.  Since $\C$ is abelian, the
restriction of $u$ to $\C$ is completely bounded on $\C$, so by the
Dixmier-Day Theorem on amenable groups, there exists an invertible
operator $S$ with $\condnum{S}\leq \norm{u|_\C}^2$ such that $(\Ad
S)\circ u|_\C$ is a $*$-homomorphism.  Since $u$ is one-to-one, we
have
$$\norm{x}^2=\norm{x^*x}=\norm{Su(x^*x)S^{-1}}\leq
\norm{u}^2\norm{u(x^*)u(x)}\leq \norm{u}^3\norm{x}\norm{u(x)},$$ so
that $\norm{x}\leq \norm{u}^3\norm{u(x)}$.  Thus, $u$ is bounded
below, so that the range of $u$ is closed.   
\end{proof}

Combining Theorems~\ref{PaulsenSimThm} and~\ref{isocor} with the
structure of completely contractive representations we obtain the
following corollary.

\begin{corollary}\label{invstructure}
Suppose $\A\subseteq \bh$ is a \cstaralg\ and $u:\A\rightarrow
\B(\H_u)$ is a faithful bounded representation.  Then there exists a
Hilbert space $\K$, a $*$-representation $\pi: \B(\H_u)\rightarrow \B(\K)$,
an isometry $W:\H\rightarrow \K$ and an invertible operator $S\in \bh$
with $\condnum{S}\leq \norm{ u^{-1}}\norm{ u}^4$ 
such that for every $x\in \A$,
\begin{equation}\label{formulaforx}
 x= SW^*\pi(u(x))WS^{-1}.
\end{equation} 

\end{corollary}

\begin{proof}
Let $\B=u(\A)$.  We may assume that $\A$ is unital and $u(I)=I$.
 Theorem~\ref{isocor} shows that $u^{-1}$ is a completely bounded map
 from $\B$ to $\bh$.  Therefore, by Theorem~\ref{PaulsenSimThm}, there
 exists an invertible operator $S\in\bh$ with
 $\condnum{S}=\norm{u^{-1}}_{\rm cb}$ so that $\psi:=\Ad (S^{-1})\circ
 u^{-1}$ is completely contractive.  By Arveson's Structure Theorem
 for completely contractive representations of operator algebras,
 there exists a Hilbert space $\K$, a $*$-representation
 $\pi:\B(\H_u)\rightarrow \B(\K)$ and an isometry $W$ so that for all
 $b\in\B$, $\psi(b)=W^*\pi(b) W$.  Letting $b=u(x)$ (for $x\in \A$)
 yields \eqref{formulaforx}.  The estimate for the condition number of
 $S$ follows from Theorem~\ref{autocb}.
\end{proof}

\begin{remark}{Remark}  Unfortunately, we have been unable to solve
  for $u(x)$ in \eqref{formulaforx}; doing so would of course lead to
  a solution of Kadison's problem.

 The range of the isometry $W$ appearing in
 Corollary~\ref{formulaforx} is a semi-invariant subspace for
 $\pi(\B)$.  Thus, $WW^*=PQ^\perp$ for some projections $P,
 Q\in\lat(\pi(\B))$, with $Q\leq P$.  The map $x\mapsto WS^{-1}x SW^*$ is a
 homomorphism into the ``$2,2$-diagonal piece'' of $\pi(u(x))$ relative to the
 block decomposition of $\pi(u(x))$ according to $I= Q+PQ^\perp +
 Q^\perp P^\perp.$
\end{remark}

We now show that Kadison problem is equivalent to the
issue of whether the image of a \cstaralg\ under a bounded
homomorphism $\lambda$ norms itself.

\begin{theorem}\label{KadEquiv}
Suppose $\A$ is a \cstaralg\ and $u:\A\rightarrow \bh$ is a bounded
representation with image $\B:=u(\A)$, and let $\tilde{u}:\A/\ker
u\rightarrow \B$ be the induced map.  If $\B$ $\lambda$-norms itself
for some $\lambda >0$, then $u$ is completely bounded and
$\norm{u}_{cb}\leq \lambda^{-1}\norm{u}^9\norm{\tilde{u}^{-1}}^2$.
Conversely, if $u$ is completely bounded, then $\B$ $\lambda$-norms
itself for any $\lambda$ with $\ds 0 < \lambda \leq \frac{1}{\norm{u}_{cb}
\norm{\tilde{u}^{-1}}\norm{u}^4}$.
\end{theorem}
\begin{proof}
We may assume that $u$ is a monomorphism, so Theorem~\ref{isocstar}
gives $\norm{u^{-1}}_{cb}\leq \norm{u}^4\norm{u^{-1}}.$  Suppose that $\B$ is
$\lambda$-normed by itself for some $\lambda>0$, and let $T\in
M_n(\A)$.  The calculation in the proof of Theorem~\ref{autocb}  shows
that if $R\in M_{1,n}(\B)$ and $C\in M_{n,1}(\B)$ satisfy
$\norm{R}\leq 1$ and $\norm{C}\leq 1$, then
$$\norm{Ru_n(T)C}\leq
\norm{u}\norm{u_{1,n}^{-1}(R)}\norm{T}\norm{u_{n,1}^{-1}(C)}\leq
\norm{u}^9\norm{u^{-1}}^2\norm{T}.$$ Taking suprema over all such $R$
and $C$ gives $$\lambda\norm{u_n(T)}\leq
\norm{u}^9\norm{u^{-1}}^2\norm{T},$$ so $u$ is completely bounded with
$\norm{u}_{cb} \leq \lambda^{-1}\norm{u}^9\norm{u^{-1}}^2.$

Conversely, suppose that $u$ is completely bounded, and view $\B$
as a bimodule over itself.  
For any $R\in \M_{1,n}(\A)$, $C\in \M_{n,1}(\A))$ with $\norm{R},
\norm{C}\leq 1,$ we
have,
\begin{align*}
\norm{Ru_n^{-1}(T)C}&\leq \norm{u^{-1}}\norm{u(Ru_n^{-1}(T)C)}\\
 &=\norm{u^{-1}}\norm{u_{1,n}(R)} \norm{u_{n,1}(C)}
 \norm{\frac{u_{1,n}(R)}{\norm{u_{1,n}(R)}} T
 \frac{u_{n,1}(C)}{\norm{u_{n,1}(C)}}}\\ &\leq \norm{u^{-1}}\norm{u}^4
 \eta_n(T) \quad \text{(using Lemma~\ref{PisierLemma}).}
\end{align*}
Taking suprema yields $\norm{u_n^{-1}(T)}\leq \norm{u^{-1}}\norm{u}^4
\eta_n(T)$, and hence
$$\norm{T}\leq \norm{u_n}\norm{u_n^{-1}(T)}\leq \norm{u}_{cb}
\norm{u^{-1}}\norm{u}^4 \eta_n(T).$$ Thus $\B$ $\lambda$-norms itself
for any $\ds \lambda \leq \frac{1}{\norm{u}_{cb}
\norm{u^{-1}}\norm{u}^4}$.
\end{proof}

\begin{corollary}\label{simC*alg}  Suppose $\B\subseteq \bh$
 is an operator algebra which is isomorphic to a \cstaralg\ $\A$.  Then
 there exists an invertible operator $S\in\bh$ such that $S\B S^{-1}$
 is a \cstaralg\ if and only if $\B$ $\lambda$-norms itself for some
 $\lambda >0$.
\end{corollary}
\begin{proof}
Let $u:\A\rightarrow \B$ be the isomorphism.  Then
Theorem~\ref{KadEquiv} shows that if $\B$ $\lambda$-norms itself, then
$u$ is completely bounded, hence by Theorem~\ref{HaagerupSimThm}, $u$
is similar to a $*$-representation, and so $\B$ is similar to a
\cstaralg.  Conversely, if $\B$ is similar to a \cstaralg, then $u$ is
similar to a $*$-representation and an application of
Theorem~\ref{KadEquiv} completes the proof.
\end{proof}

The following question is thus a reformulation of Kadison's question.

\begin{remark}{Question}  Suppose $\B$ is an operator algebra which is
  isomorphic to a \cstaralg.  Does $\B$ $\lambda$-norm itself for some
  $\lambda>0$?
\end{remark}

\begin{remark}{Remark}  Haagerup~\cite{HaagerupSoSiPrCyRe}
  showed that every bounded, cyclic representation $u$ of a \cstaralg\
  $\A$ is completely bounded with $\norm{u}_{cb}\leq \norm{u}^3$.  So
  given an arbitrary representation $u$ of $\A$ on $\H$, let
  $\B=u(\A)$.  For each unit vector $\xi\in\H$, let $P_\xi$ be the
  projection onto the cyclic subspace $[\B\xi]$ and let $\B_\xi$ be
  the restriction of $\B$ to $\H_\xi:=P_\xi\H$. The representation
  $u_\xi$ of $\A$ given by $u_\xi(x)=u(x)P_\xi$ is thus completely
  bounded and $\norm{u_\xi}_{cb}\leq \norm{u}_{cb}$.
  Theorem~\ref{PaulsenSimThm} shows that there exists an invertible
  operator $S_\xi\in\B(\H_\xi)$ such that $\norm{u_\xi}_{cb}=
  \condnum{S_\xi}$ and $(\Ad S_\xi) (\B_\xi)$ is a \cstaralg, call it
  $\A_\xi$.  Applying Theorem~\ref{KadEquiv} to $\Ad
  S^{-1}_\xi:\A_\xi\rightarrow \B_\xi$, we find $\norm{\Ad S_\xi}\leq
  \condnum{S_\xi}$, so that $\B_\xi$ $\lambda$-norms itself with
  $\lambda= \norm{u}^{-18}.$
\end{remark}

In the following example, we show that there exists an operator
algebra which does not $\lambda$-norm itself.  The idea for the proof
is due to Ken Davidson.

\begin{remark}{Example}  Let $\bbD\subseteq \bbC$ be the open unit 
  disk and let $\bbA(\bbD)\subseteq C(\overline{\bbD})$ be the disk
  algebra, that is, the collection of all continuous functions on the
  closed unit disk which are analytic in $\bbD$.  We use $\P\subseteq
  \bbA(\bbD)$ to denote the collection of all polynomials.

  Recall that an operator $T\in\bh$ is polynomially bounded if there
  exists $K>0$ such that for every $p\in\P$, $\norm{p(T)}\leq
  K\norm{p}$.  Polynomial boundedness of $T$ is equivalent to the
  existence of a bounded homomorphism $u:\bbA(\bbD)\rightarrow \bh$
  such that for every $p\in\P,$ $u(p)=p(T)$.  A polynomially bounded
  operator $T$ is \textit{completely polynomially bounded} if $u$ is
  completely bounded.  Paulsen~\cite{PaulsenEvCoPoBoOpSiCo} showed
  that $T$ is completely polynomially bounded if and only if $T$ is
  similar to a contraction.  Pisier~\cite{PisierAPoBoHiSpNoSiCo}
  showed that there exists a polynomially bounded operator $T\in\bh$
  which is not completely polynomially bounded, so $T$ is not similar
  to a contraction.

  Fix a polynomially bounded operator $T\in\bh$.  Notice that the
  spectrum of $T$ is contained in $\overline{\bbD}$.  If $U$ is
  a unitary operator with $\sigma(U)=\bbT$, then $T$
  is completely polynomially bounded if and only if $T\oplus U$ is
  completely polynomially bounded.   We may therefore assume that
  $$\bbT\subseteq \sigma(T).$$ Let $u:\bbA(\bbD)\rightarrow\bh$ be the
  extension of the map $p\in\P\mapsto p(T)$ and set
  $\B=u(\bbA(\bbD))$.  Since $\bbT\subseteq \sigma(T)$, we have
  $\norm{u(f)}\geq \norm{f}$ for every $f\in\bbA(\bbD)$, so that $\B$
  is closed and $u^{-1}$ is contractive.  Since $\bbA(\bbD)\subseteq
  C(\overline{\bbD})$, the operator space structure on $\bbA(\bbD)$ is
  minimal among all all possible operator space structures on
  $\bbA(\bbD)$ (see~\cite[Paragraph~1.2.21]{BlecherLeMerdyOpAlThMo}),
  hence $u^{-1}$ is completely contractive.

  View $\B$ as a $\B$-bimodule and let $\eta_n$ be the norm on $M_n(\B)$ as in
  Definition~\ref{normingdef}.  We shall show that for every $n$,
  the norm $\eta_n(u_n(\cdot))$ and the usual norm
  $\norm{\cdot}_{M_n(\bbA(\bbD))}$ are equivalent norms on
  $M_n(\bbA(\bbD))$.  

  Choose $X\in M_n(\bbA(\bbD))$, $R\in M_{1n}(\bbA(\bbD))$ and $C\in
  M_{n1}(\bbA(\bbD))$.  Since $u$ is bounded, we have
  $\norm{u_{1n}(R)u_n(X)u_{n1}(C)}=\norm{u(RXC)}\leq
  \norm{u}\norm{RXC}$.  Since $\norm{u_{1n}(R)}\geq \norm{R}$ and
  $\norm{u_{n1}(C)}\geq \norm{C}$, we obtain
  \begin{equation}\label{first}
  \eta_n(u_n(X))\leq \norm{u} \norm{X}.
 \end{equation}
  On the other hand, a result of
  Bourgain~\cite{BourgainNeBaSpPrDiAlHIn} (see
  also~\cite[Theorem~9.9]{PisierSiPrCoBoMaIIed}) shows that there
  exists a constant $s>0$ (independent of $u$) such that  for every
  $n$, $\max\{\norm{u_{1n}}, \norm{u_{n1}}\}\leq s\norm{u}$.  Therefore,
\begin{align*}
  \norm{RXC}&\leq \norm{u_{1n}(R)u_n(X)u_{n1}(C)}\\
   &\leq s^2\norm{u}^2\norm{\frac{u_{1n}(R)}{\norm{u_{1n}(R)}}
   u_n(X)\frac{u_{n1}(C)}{\norm{u_{n1}(C)}}}\norm{R}\norm{C}.
\end{align*}  
Since the operator space structure on $\bbA(\bbD)$ is minimal,
$\bbA(\bbD)$ norms itself.  Taking the supremum over $R$ and $C$ with
norm one, we obtain 
 \begin{equation}\label{second}
\norm{X}\leq s^2\norm{u}^2\eta_n(u_n(X)).
\end{equation}
Combining \eqref{first} and \eqref{second} establishes the claim.

Thus, when $T$ is chosen to be polynomially bounded but not completely
polynomially bounded, we see that $\B$ cannot $\lambda$-norm itself.
\end{remark}

\subsection{An Application to CSL Algebras}
Recall that a \textit{CSL algebra} is an operator algebra $\A\subseteq \bh$
which is both reflexive and such that there exists a MASA $\D\subseteq
\bh$ with $\D\subseteq \A$.  The lattice of invariant projections of a
CSL algebra is a commutative family of projections, and when this
lattice is completely distributive, the CSL algebra is called a
\textit{completely distributive} CSL algebra.

Davidson and Power proved that isometric isomorphisms of completely
distributive CSL algebras are unitarily implemented.  Their techniques
involved   homological ideas and were
somewhat intricate.  We can give a simpler proof of their result,
and which also extends theirs.  For convenience, we assume
irreducibilty.  If $\A$ is not irreducible, one can use 
direct integrals along the center of $\A\cap \A^*$, together with an
appropriate hypothesis on the lattices of each term in the direct
integral to obtain a more general result.

\begin{theorem}\label{isomCSL}  For $i=1,2$, let $\A_i
  \subseteq\B(\H_i)$ be  CSL
  algebras such that $\A_i\cap\K_i\neq (0)$  and suppose that $\A_1$
  is irreducible.   If $ u:\A_1\rightarrow \A_2$ is an isometric isomorphism,
  then there exists a unitary operator $U\in \B(\H_1,\H_2)$ such that
  $ u=\Ad U$.
\end{theorem}

\begin{proof}
Since the isomorphism $ u$ is isometric, its restriction to
$\A_1\cap (\A_1)^*$ is a $*$-isomorphism, so in particular, $\lat(\A_2)$ is
isomorphic to $\lat(\A_1)$.   Thus, since $\A_1$ is irreducible, so is
$\A_2$.  Let $\C_i$ be the \cstar-subalgebra of $\B(\H_i)$ generated
by $\A_i$.  Then $\C_i$ is an irreducible \cstaralg\ containing a
compact operator, hence $\C_i$ contains all the compact operators.
There exists a $*$-epimomorphism $\pi:\C_i\rightarrow C^*_{\rm
env}(\A_i)$ such that $\pi|_{\A_i}=\iota$, where $\iota$ is the
canonical inclusion of $\A_i$ into $C^*_{\rm env}(\A_i)$.  If
$\ker\pi\neq (0)$, then since $\C_i$ contains the compact operators,
$\ker\pi\cap\K_i\neq (0).$ Hence $\K_i\subseteq \ker\pi$, since
$\ker\pi\cap \K_i$ is an ideal in an irreducible \cstaralg.  But
this is impossible, since $\pi$ is isometric on $\A_i\cap \K_i$.
Therefore $\ker\pi=(0)$, so that $\C_i$ is the \cstar-envelope of
$\A_i$.

Theorem~2.7 of~\cite{PopSinclairSmithNoC*Al} shows that any MASA is
norming for $\bh$, hence $\A_i$ contain norming \cstar-subalgebras.  
Theorem~\ref{autocb} shows that $ u$ and $ u^{-1}$ are
completely contractive, so that $ u$ is a complete isometry.  By
the universal property of \cstar-envelopes (applied to $ u$ and
$ u^{-1}$), $ u$ extends to a $*$-isomorphism $\tilde{ u}$
of $\C_1$ onto $\C_2$.  The compact operators
are the smallest closed two-sided ideal contained in
$\C_i$, so that the restriction of $\tilde{ u}$ to the
compact operators is a $*$-isomorphism of $\K(\H_1)$ onto $\K(\H_2)$.
Therefore, there exists a unitary operator $U$ so that $(\Ad
U)|_{\K(\H_1)}=\tilde{ u}|_{\K(\H_1)}$.  Finally, if
$T\in\C_1$ and if $\eta\in\H_2$, we may find a finite rank
projection $P$ so that
$2\norm{T}\norm{P^\perp \eta} <\eps$.
Then since $\tilde{ u}^{-1}(P)=(\Ad U^*) (P)$, we have
\begin{align*}
\norm{\left(\tilde{ u}(T)-(\Ad U)(T)\right)\eta} &\leq
\norm{\left(\tilde{ u}(T\tilde{ u}^{-1}(P))- 
(\Ad U)(T (\Ad U^*(P))\right)\eta} \\ &\qquad + 
\norm{\left(\tilde{ u}(T)-(\Ad U)(T)\right)P^\perp\eta} \\
&= \norm{\left(\tilde{ u}(T)-(\Ad U)(T)\right)P^\perp\eta}<\eps,
\end{align*} so
$\tilde{ u}(T) = (\Ad U)(T).$  
Since $\A_1\subseteq \C_1$,
the proof is complete.
\end{proof}

\subsection{Applications to Subalgebras of \cstar-Diagonals}
In this subsection, we provide applications to subalgebras of certain
classes of \cstaralg s.  

A \cstar-diagonal is a pair $(\C,\D)$ of \cstaralg s such that $\D$ is
abelian and such that  
\begin{enumerate}
\item[i)] every pure
state of $\D$ extends uniquely to a pure state of $\C$;
\item[ii)] the conditional expectation $E:\C\rightarrow \D$ (whose
  existence is guaranteed by (i)) is faithful;
\item[iii)] the closed linear span of the set $\{v\in\C: v\D=\D v\}$
  is $\C$.
\end{enumerate}  We will assume that both $\C$ and $\D$ are unital.
  The extension property then 
  implies that $\D$ is a MASA in $\C$.  

Such pairs were introduced by
  Kumjian~\cite{KumjianOnC*Di}, who used slightly different, but
  essentially equivalent axioms (see~\cite{DonsigPittsCoSyBoIs} for a
  discussion of the equivalence).  Also, \cstar-diagonals and their subalgebras
  were further in several papers, see for
  example~\cite{DonsigPittsCoSyBoIs,MuhlyQiuSolelCoNuSpSuOpAl,%
  MuhlyQiuSolelIsOpAl}.  

Our first task is to show that $\D$ norms $\C$.  
\begin{lemma}\label{dnorming}
Suppose $(\C,\D)$ is a \cstar-diagonal.  Then $\C$ is normed by $\D$.
\end{lemma}
\begin{proof}
Theorem~5.9 of~\cite{DonsigPittsCoSyBoIs} shows that there exists a
faithful $*$-representation $\pi:\C\rightarrow\bh$ such that
$\pi(\D)''$ is a MASA in \bh.  It follows from~\cite[Lemma~2.2 and
  Theorem~2.7]{PopSinclairSmithNoC*Al} that $\pi(\D)$ norms \bh, hence
$\pi(\D)$ norms $\pi(\C)$.  As $\pi$ is a faithful $*$-representation of a
\cstaralg, it is a complete
isometry, so $\D$ norms $\C$.
\end{proof}  

The following notation will be useful.  When $(\C,\D)$ is a
\cstar-diagonal and $\A$ is a norm closed algebra with $\D\subseteq
\A\subseteq \C$, we will write $\A\subseteq (\C,\D)$.

For $i=1,2$, let $(\C_i,\D_i)$ be \cstar-diagonals.  Muhly, Qiu and
Solel~\cite[Theorem~1.1]{MuhlyQiuSolelIsOpAl}  proved that
when $\A_i\subseteq (\C_i,\D_i)$ are triangular, that is $\A_i\cap
(\A_i)^*=\D_i$,  which generate $\C_i$
and $(\C_i,\D_i)$ are nuclear, then an isometric isomorphism
$ u:\A_1\rightarrow \A_2$ extends to a $*$-isomorphism of $\C_1$
onto $\C_2$.  Later Donsig and
Pitts~\cite[Theorem~8.9]{DonsigPittsCoSyBoIs} extended this result:
they showed that the hypothesis of nuclearity can be removed.
To prove their result, Donsig and
Pitts used showed that the isometric isomorphism between $\A_1$ and
$\A_2$ induces 
isomorphism of an appropriate CSL
algebras, then used the structure theory for isomorphisms of CSL
algebras.
The techniques used to prove
~\cite[Theorem~1.1]{MuhlyQiuSolelIsOpAl}
and~\cite[Theorem~8.9]{DonsigPittsCoSyBoIs} do not apply for
non-triangular subalgebras.  A

Donsig and Pitts~\cite[Theorem~4.22]{DonsigPittsCoSyBoIs} showed that
the \cstar-envelope of any subalgebra (triangular or not) $\A\subseteq
(\C,\D)$ is the \cstar-subalgebra, $C^*(\A)$, of $\C$ generated by
$\A$.  In the context of both~\cite[Theorem~1.1]{MuhlyQiuSolelIsOpAl}
and~\cite[Theorem~8.9]{DonsigPittsCoSyBoIs}, $(C^*(\A_i),\D_i)$ are
\cstar-diagonals.  When $\C$ is separable and nuclear, the Spectral
Theorem for Bimodules~\cite{MuhlyQiuSolelCoNuSpSuOpAl} shows that
$(C^*(\A),\D)$ is a \cstar-diagonal.  In the general case however, it
is not clear that the pair $(C^*(\A),\D)$ is a \cstar-diagonal---one
needs to verify that condition (ii) of the definition of
\cstar-diagonal holds.  The following consequence of
\cite[Theorem~4.22]{DonsigPittsCoSyBoIs}, Corollary~\ref{isometriciso}
and Lemma~\ref{dnorming} is therefore a significant extension of
\cite[Theorem~1.1]{MuhlyQiuSolelIsOpAl}
and~\cite[Theorem~8.9]{DonsigPittsCoSyBoIs}.

\begin{theorem}\label{cstardiag}
Let $\A_i\subseteq (\C_i,\D_i)$ be norm-closed subalgebras of
\cstar-diagonals.  If $ u:\A_1\rightarrow\A_2$ is an isometric
isomorphism, then $ u$ extends uniquely to a $*$-isomorphism of
$C^*(\A_1)$ onto $C^*(\A_2)$.  
\end{theorem}

\begin{remark}{Remark}  In~\cite{MercerIsIsCaBiAl}, Mercer 
  proves a result similar to Theorem~\ref{cstardiag}, but where the
  algebras $\A_i$ are taken to be weak-$*$ closed subalgebras of von
  Neumann algebras $\fM_i$ and there are Cartan MASAs $\D_i\subseteq
  \fM_i$ such that $\D_i\subseteq \A_i\subset\M_i$.  We expect that
  Cartan MASAs norm their containing von Neumann algebras, and thus expect
  that it should be possible to give a proof of Mercer's result based
  on Theorem~\ref{autocb} as well.
\end{remark}

\vskip .5em \noindent\textsc{Acknowledgment:} This work was completed during a
visit to the University of Waterloo.  The author thanks Ken Davidson
and the University of Waterloo for their hospitality, and also thanks
Ken Davidson and Allan Donsig for several helpful conversations and
suggestions.

\bibliographystyle{amsplain}


\providecommand{\bysame}{\leavevmode\hbox to3em{\hrulefill}\thinspace}
\providecommand{\MR}{\relax\ifhmode\unskip\space\fi }
\newcommand{\mr}{\relax}
\providecommand{\MRhref}[2]{%
  \href{http://www.ams.org/mathscinet-getitem?mr=#1}{#2}
}
\providecommand{\href}[2]{#2}

\end{document}